\documentclass{anmatuvt} 
\usepackage{mathrsfs, amssymb, eucal, amsmath, amsthm, amscd, amsrefs}
\numberwithin{equation}{section}
\newtheorem{thm}{\bf Theorem}[section]
\newtheorem{lem}[thm]{\bf Lemma}
\newtheorem{prop}[thm]{\bf Proposition}
\newtheorem{cor}[thm]{\bf Corollary}
\newtheorem{defn}{\bf Definition}[section]
\theoremstyle{remark}
\newtheorem{rem}{\bf Remark}[section]
\newtheorem{exmp}{\bf Example}[section]

\begin{document}
\title[$\eta$-Ricci Solitons on Kenmotsu 3-Manifolds]{$\eta$-Ricci Solitons on Kenmotsu 3-Manifolds}

\author[Krishnendu De]{Krishnendu De }
\address{
Assistant Professor of Mathematics, \\
Kabi Sukanta Mahavidyalaya, \\
Bhadreswar, P.O.-Angus, Hooghly, \\
Pin 712221, West Bengal,\\
 India.}
\email{krishnendu.de@outlook.in}

\author[Uday Chand De]{Uday Chand De}
\address{
Department of Pure Mathematics, \\
University of
Calcutta, \\
35, Ballygaunge Circular Road, \\
Kolkata 700019, \\
India}
\email{uc$\_$de@yahoo.com}

\date{ }
\thanks{}
\dedicatory{}

\maketitle

\subjclass{53C15 ; 53C25.}
\keywords{Ricci soliton, $\eta$-Ricci soliton, Kenmotsu 3-manifolds, Coddazi type of Ricci tensor, cyclic parallel Ricci tensor.}

\begin{abstract}
In the present paper we study $\eta$-Ricci solitons on Kenmotsu 3-manifolds. Moreover, we consider $\eta$-Ricci solitons on Kenmotsu 3-manifolds with Codazzi type of Ricci tensor and cyclic parallel Ricci tensor. Beside these, we study $\phi$-Ricci symmetric $\eta$-Ricci soliton on Kenmotsu 3-manifolds. Also Kenmotsu 3-manifolds satisfying the curvature condition $R.R=Q(S,R)$ is considered. Finally, an example is constructed to prove the existence of a proper $\eta$-Ricci soliton on a Kenmotsu 3-manifold.
\end{abstract}

\section{Introduction}
\hspace{.3cm} In modern mathematics, the methods of contact geometry play an important role. Contact geometry has evolved from the mathematical formalism of classical mechanics. In the present paper we are entering an era of new concepts, some generalizations and applications which play a functional role in contemporary mathematics. The product of an almost contact manifold $M$ and the real line $\mathbb {R}$ carries a natural almost complex structure. However, if one takes $M$ to be an almost contact metric manifold and suppose that the product metric $G$ on $M\times \mathbb {R}$ is Kaehlerian, then the structure on $M$ is cosymplectic \cite{ian} and not Sasakian. On the other hand, Oubina \cite{oub} pointed out that if the conformally related metric $e^{2t}G$, $t$ being the coordinate on $\mathbb {R}$, is Kaehlerian, then $M$ is Sasakian and conversely.\par

\hspace{.3cm}In \cite{tan} S. Tanno classified almost contact metric manifolds whose automorphism groups possess the maximum dimension. For such a manifold $M$, the sectional curvature of plane section containing $\xi$ is a constant, say $c$. If $c>0$, $M$ is a homogeneous Sasakian manifold of constant sectional curvature. If $c=0$, $M$ is the product of a line or a circle with a Kaehler manifold of constant holomorphic sectional curvature. If $c<0$, $M$ is a warped product space $\mathbb {R}\times_{ f}C ^{n}$. In $1972$, K. Kenmotsu \cite{ken} abstracted the differential geometric properties of the third case. We call it Kenmotsu manifold.\par
\hspace{.3cm} In 1982, R. S. Hamilton \cite{hr} introduced the notion of Ricci flow to find a canonical metric on a smooth manifold. The Ricci flow is an evolution equation for metrics on
a Riemannian manifold defined as follows:
\begin{eqnarray} \frac{\partial}{\partial_t}g_{ij}=-2R_{ij}.\label{a1}\end{eqnarray}
Ricci solitons are special solutions of the Ricci flow equation (\ref{a1}) of the form $g_{ij}=\sigma(t)\psi_t^\ast g_{ij}$ with the initial condition $g_{ij}(0)=g_{ij}$, where $\psi_t$ are diffeomorphisms of $M$ and $\sigma(t)$ is the scaling function.\par

\hspace{.3cm} A Ricci soliton is a generalization of an Einstein metric. We recall the notion of Ricci soliton according to \cite{CC}. On the manifold  $M$, a Ricci soliton is a triple $(g,V,\lambda )$ with $g$ a Riemannian metric, $V$ a vector field, called potential vector field  and $\lambda $ a real scalar such that
 \begin{eqnarray} \pounds_{V}g+2S+2\lambda g=0,\label{a2}\end{eqnarray}
where $\pounds$ is the Lie derivative. Metrics satisfying (\ref{a2}) are interesting and useful in physics and are often referred as quasi-Einstein (\cite{chave1},\cite{chave2}).\par

\hspace{.3cm} Theoretical physicists have also been looking into the equation of Ricci soliton in relation with string theory. The initial
contribution in this direction is due to Friedan \cite{friedam} who discusses some aspects of it.\par

\hspace{.3cm} The Ricci soliton is said to be shrinking, steady and expanding according as $\lambda $ is negative, zero and positive respectively. Ricci solitons have been studied by several authors such
as (\cite{deshmukh1},\cite{deshmukh2},\cite{hr},\cite{dugg},\cite{it}) and many others.\par

\hspace{.3cm} As a generalization of Ricci soliton, the notion of $\eta$-Ricci soliton was introduced by Cho and Kimura \cite{cho1}. This notion has also been studied in \cite{CC} for Hopf hypersurfaces in complex space forms. An $\eta$-Ricci soliton is a tuple $(g,V,\lambda ,\mu )$, where $V$ is a vector field on $M$, $\lambda$ and $\mu$ are constants, and $g$ is a Riemannian (or pseudo-Riemannian) metric satisfying the equation
\begin{eqnarray}
   \pounds_{V}g+2S+2\lambda g+2\mu \eta \otimes \eta =0,\label{a3}
  \end{eqnarray}
where $S$ is the Ricci tensor associated to $g$. In this connection, we mention the works of Blaga (\cite{blaga1}, \cite{blaga2}) and Prakasha et al. \cite{dg} on $\eta$-Ricci solitons. In particular, if $\mu =0$, then the notion of $\eta$-Ricci soliton $(g,V,\lambda ,\mu )$ reduces to the notion of Ricci soliton $(g,V,\lambda )$. If $\mu\not=0$, then the $\eta$-Ricci soliton is named proper $\eta$-Ricci soliton.\\

Motivated by the above works, we study $\eta$-Ricci solitons on Kenmotsu 3-manifolds.\\

The paper is organized as follows:\par
After preliminaries in Section 2, we study $\eta$-Ricci solitons on a Kenmotsu 3-manifold. Moreover, we consider $\eta$-Ricci solitons on Kenmotsu 3-manifolds with Codazzi type of Ricci tensor and cyclic parallel Ricci tensor. Beside these, we study $\phi$-Ricci symmetric $\eta$-Ricci solitons on Kenmotsu 3-manifolds. In the next section, Kenmotsu 3-manifolds satisfying the curvature condition $R.R=Q(S,R)$ are studied. Finally, an example is constructed to prove the existence of a proper $\eta$-Ricci soliton on a Kenmotsu 3-manifold.

\section{Preliminaries}

\hspace{.3cm} Let $M$ be a connected almost contact metric manifold with an almost contact metric structure ($\phi$,$\xi$,$\eta$,$g$), that is, $\phi$ is an ($1$,$1$)-tensor field, $\xi$ is a vector field, $\eta$ is a $1$-form and $g$ is a compatible Riemannian metric such that
\begin{eqnarray}\label{b1}
\phi^{2}(X)=-X +\eta(X)\xi, \eta(\xi)= 1, \phi\xi= 0, \eta\phi  =
0
\end{eqnarray}
\begin{eqnarray}\label{b2}
 g(\phi X,\phi Y) =
g(X,Y) - \eta(X)\eta(Y)
\end{eqnarray}
\begin{eqnarray}\label{b3}
g(X,\xi) = \eta(X)
\end{eqnarray}
for all $X,Y\in \Gamma (TM)$(\cite{bla1},\cite{bla2}).\\
\hspace{.3cm} If an almost contact metric manifold satisfies
\begin{eqnarray} (\nabla _{X}\phi)Y=g(\phi X,Y)\xi-\eta (Y)\phi X,\label{b4}\end{eqnarray}
then $M$ is called a Kenmotsu manifold \cite{ken}, where $\nabla$ is the Levi-Civita connection of $g$. From the above equation it follows that
\begin{eqnarray} \nabla _{X}\xi=X- \eta (X)\xi\label{b5}\end{eqnarray} and
\begin{eqnarray}(\nabla _{X}\eta)Y= g(X,Y)-\eta (X)\eta (Y).\label{b6}\end{eqnarray}
\hspace{.3cm} Moreover, the curvature tensor $R$ and the Ricci tensor $S$ satisfy
\begin{eqnarray} R(X,Y)\xi= \eta (X)Y- \eta (Y)X,\label{b7}\end{eqnarray}
\begin{eqnarray} R(\xi,X)Y= \eta (Y)X-g(X,Y)\xi,\label{b8}\end{eqnarray}
\begin{eqnarray} R(\xi,X)\xi=X-\eta (X)\xi,\label{b9}\end{eqnarray}
and \begin{eqnarray} S(X,\xi)= -(n-1)\eta (X).\label{b10}\end{eqnarray}
\hspace{.3cm} From \cite{de}, we know that for a 3-dimensional Kenmotsu manifold
\begin{eqnarray} R(X,Y)Z &=&\frac{r+4}{2}[g(Y,Z)X-g(X,Z)Y]\\&&
-\frac{r+6}{2}[g(Y,Z)\eta (X)\xi-g(X,Z)\eta (Y)\xi\nonumber\\&&+\eta (Y)\eta (Z)X-\eta (X)\eta
(Z)Y]\nonumber,\label{b11}\end{eqnarray}
\begin{eqnarray} S(X,Y)=\frac{1}{2}[(r+2)g(X,Y)-(r+6)\eta (X)\eta (Y)],\label{b12}\end{eqnarray}
where $S$ is the Ricci tensor of type (0,2), $R$ is the curvature tensor of type (1,3) and $r$ is the scalar curvature of the manifold $M$.

{\prop For an $\eta$-Ricci soliton on a Kenmotsu 3-manifold we have $\lambda +\mu =2$.}
\proof Assume that the Kenmotsu $3$-manifold admits a proper $\eta$-Ricci soliton $(g,\xi ,\lambda ,\mu )$. Then the relation $(1.3)$ yields
\begin{eqnarray}
(\pounds_{\xi }g)(X,Y)+2S(X,Y)+2\lambda g(X,Y)+2\mu \eta (X)\eta (Y)=0 .\label{b13}\end{eqnarray}
It follows that
\begin{eqnarray}
 2S(X,Y)=-(\pounds _{\xi }g)(X,Y)-2\lambda g(X,Y)-2\mu \eta (X)\eta (Y),\label{b14}\end{eqnarray}
for all smooth vector fields $X,Y\in \Gamma (TM)$. In a Kenmotsu 3-manifold we have
\begin{eqnarray} (\pounds _{\xi }g)(X,Y)=2[g(X,Y)-\eta(X)\eta(Y)].\label{b15}\end{eqnarray}
Making use of (\ref{b15}) in (\ref{b14}) we get
\begin{eqnarray} S(X,Y)=-(\lambda+1) g(X,Y)-(\mu-1) \eta (X)\eta (Y)\label{b16}.\end{eqnarray}

Comparing the above equation (\ref{b12}) with (\ref{b16}), we find $\lambda+1 =-\frac{1}{2}(r+2)$ and $\mu-1 =\frac{1}{2}(6+r),$ from which it follows that $\lambda +\mu =2.$ This completes the proof.

\section{$\eta$-Ricci Solitons on Kenmotsu 3-manifolds with Coddazi type of Ricci tensor}
In this section we consider proper $\eta$-Ricci soliton on Kenmotsu 3-manifolds with Coddazi type of Ricci tensor. A. Gray \cite{gray} introduced the notion of cyclic parallel Ricci tensor and Codazzi type of Ricci tensor.\\

\hspace{.3cm} A Riemannian manifold is said to satisfy cyclic parallel Ricci tensor if its Ricci tensor $S$ of type (0,2) is non-zero and satisfies  the condition
\begin{eqnarray}
 (\nabla _{X}S)(Y,Z)+(\nabla _{Y}S)(Z,X)+(\nabla _{Z}S)(X,Y)=0.\label{c1}\end{eqnarray}

\hspace{.3cm} A Riemannian manifold is said to have Codazzi type of Ricci tensor if its Ricci tensor $S$ of type (0,2) is non-zero and satisfies the following condition
\begin{eqnarray} (\nabla _{X}S)(Y,Z)=(\nabla _{Y}S)(X,Z).\label{c2}\end{eqnarray}

Therefore, taking covariant differentiation of (\ref{b16}) with respect to $Z$ we obtain
\begin{eqnarray}
(\nabla_{Z}S)(X,Y)&=&-(\mu-1) [(\nabla_{Z}\eta )(X)\eta (Y)+\eta (X)(\nabla_{Z}\eta )(Y)] \nonumber\\
  &=&-(\mu-1)[g(X,Z)\eta (Y)+g(Y,Z)\eta (X)\nonumber\\&-&2\eta(X)\eta(Y)\eta(Z)]. \label{c3}\end{eqnarray}

If the Ricci tensor $S$ is of Coddazi type, then
\begin{eqnarray} (\nabla_{Z}S)(X,Y)=(\nabla_{X}S)(Z,Y).\label{c4} \end{eqnarray}

Using (\ref{c3}) in (\ref{c4}) we have
\begin{eqnarray}
  &&(\mu-1)[g(X,Z)\eta (Y)+g(Y,Z)\eta (X)-2\eta(X)\eta(Y)\eta(Z)]\nonumber\\
  &=&(\mu-1)[g(Z,X)\eta (Y)+g(Y,X)\eta (Z)-2\eta(X)\eta(Y)\eta(Z)].\label{c5}
  \end{eqnarray}

We get $\mu=1$, hence from Proposition 2.1 we have $\lambda=1$ and from (\ref{b16}) we get $S=-2g$, from which we can easily obtain $r=-6$, where $r$ is the scalar curvature of the manifold. In a 3-dimensional Riemannian manifold, we have
\begin{eqnarray} R(X,Y)Z &=&g(Y,Z)QX-g(X,Z)QY
+S(Y,Z)X-S(X,Z)Y \nonumber\\&&-\frac{r}{2}[g(Y,Z)X-g(X,Z)Y],\label{c6}\end{eqnarray}

where $Q$ is the Ricci operator, that is $g(QX,Y)=S(X,Y)$. Now using the values of $S$, $Q$ and $r$ in the above equation we get
\begin{eqnarray} R(X,Y)Z =-[g(Y,Z)X-g(X,Z)Y].\label{c7}\end{eqnarray}

 Thus we conclude the following:
 \begin{thm} A Kenmotsu 3-manifold with Coddazi type of Ricci tensor admittting a proper $\eta$-Ricci soliton of the type $(g,V,1, 1 )$  is locally isometric to the hyperbolic space  H$(-1)$.\end{thm} 

\section{$\eta$-Ricci solitons on Kenmotsu 3-manifolds with cyclic parallel Ricci tensor}
\hspace{.3cm} This section is devoted to study proper $\eta$-Ricci solitons on Kenmotsu 3-manifolds with cyclic parallel Ricci tensor. Therefore
\begin{eqnarray}
(\nabla_{X}S)(Y,Z)+(\nabla_{Y}S)(Z,X)+(\nabla_{Z}S)(X,Y)=0,\label{d1}
\end{eqnarray}
for all smooth vector fields $X, Y, Z\in \Gamma (TM)$.
Using (\ref{c3}) in (\ref{d1}), we have
\begin{eqnarray}
 &&(\mu-1)[2g(X,Y)\eta (Z)+2g(X,Z)\eta (Y)\nonumber\\
   &&+2g(Y,Z)\eta (X)-6\eta(X)\eta(Y)\eta(Z)]=0.\label{d2}
   \end{eqnarray}
Putting $X=\xi$ in (\ref{d2}), we get
 \begin{eqnarray}
  (\mu-1) [g(Y,Z)-\eta(Y)\eta(Z)]=0. \nonumber
  \end{eqnarray}
It follows that
\begin{eqnarray} \mu=1.\label{d3}\end{eqnarray}
 Thus we, like in the earlier section, are in a position to state the following:

\begin{thm} A Kenmotsu 3-manifold with cyclic parallel Ricci tensor admittting a proper $\eta$-Ricci soliton of the type $(g,V,1, 1 )$ is locally isometric to the hyperbolic space  H$(-1)$.\end{thm}

\section{ $\phi$-Ricci Symmetric $\eta$-Ricci solitons on Kenmotsu 3-manifolds}
\hspace{.3cm} In this section we study $\phi$-Ricci symmetric proper $\eta$-Ricci solitons on Kenmotsu 3-manifolds. A Kenmotsu manifold is said to be $\phi$-Ricci symmetric if
\begin{eqnarray}
\phi^{2}(\nabla_{X}Q)Y=0,\label{e1}
 \end{eqnarray}
 for all smooth vector fields $X$, $Y$. \\
 The Ricci tensor for an $\eta$-Ricci soliton on Kenmotsu 3-manifold is given by
 \begin{eqnarray}
  S(X,Y)=-(\lambda+1) g(X,Y)-(\mu-1) \eta (X)\eta (Y).\label{e2}
  \end{eqnarray}
Then it follows that
  \begin{eqnarray}
    QX=-(\lambda+1) X-(\mu-1) \eta (X)\xi,\label{e3}
    \end{eqnarray}
for all smooth vector fields $X$.

Replacing $Q$ from (\ref{e3}) in $(\nabla_{X}Q)Y= \nabla_{X}QY -Q(\nabla_{X}Y)$ and applying $\phi^2$
we obtain:
\begin{eqnarray} (\mu-1)\eta(Y )[X -\eta(X)\xi] = 0.\nonumber\end{eqnarray}
It follows $\mu = 1$, $\lambda = 1$ and $S =-2g$, from which we can easily obtain $r=-6$.
Thus we can state the following:
\begin{thm} Let $(M,\phi, \xi,\eta,g)$ be a Kenmotsu 3-manifold. If $M$ is $\phi$-Ricci symmetric, then $\mu=1$, $\lambda=1$ and the manifold is locally isometric to the hyperbolic space  H$(-1)$.\end{thm}

 \section{$\eta$-Ricci solitons on Kenmotsu 3-manifolds satisfying $R.R=Q(S,R)$}

\hspace{.3cm} In this section we deal with $\eta$-Ricci solitons on Kenmotsu 3-manifolds satisfying $R.R=Q(S,R)$. If the tensors $R.R$ and $Q(S,R)$ are linearly dependent, then $M$ is called Ricci generalized pseudo-symmetric. This is equivalent to
\begin{eqnarray}\nonumber
R.R=fQ(S,R),
\end{eqnarray}
holding on the set $U_{R}=\{x\in M:R\neq0 $\;at\;$ x\},$ where $f$ is some function on $U_{R}.$\
A very important subclass of this class of manifolds realizing the condition is
\begin{eqnarray}\label{ee1}
R.R=Q(S,R).
\end{eqnarray}
The manifolds satisfying the condition $R.R=Q(S,R)$ were considered in (\cite{dd1},\cite{dd2}). Conformally flat manifolds realizing $(\ref{ee1})$ were investigated in \cite{dd3}. Also every 3-dimensional Riemannian manifold satisfies the above equation identically \cite{dd3}.
Now from (\ref{ee1}) we have
\begin{eqnarray}\label{ee2}
R.R=Q(S,R),
\end{eqnarray}
that is,
\begin{eqnarray}\label{ee3}
(R(X,Y).R)(U,V)W=((X\wedge _{S}Y)\cdot R)(U,V)W.
\end{eqnarray}
We get from $(\ref{ee3})$
\begin{eqnarray}\label{ee4}
&&R(X,Y)R(U,V)W-R(R(X,Y)U,V)W\nonumber\\&&-R(U,R(X,Y)V)W-R(U,V)R(X,Y)W\nonumber\\&&=(X\wedge_{S}Y)R(U,V)W-R((X\wedge_{S}Y)U,V)W
\nonumber\\&&-R(U,(X\wedge_{S}Y)V)W-R(U,V)(X\wedge_{S}Y)W.
\end{eqnarray}
We define endomorphisms $X\wedge_{A} Y$ by
\begin{eqnarray}\label{ee5}
(X\wedge_{A} Y)Z=A(Y,Z)X-A(X,Z)Y,
\end{eqnarray}
where $X,Y,Z\in\Gamma(TM)$ and $A$ is a symmetric $(0,2)$-tensor.
In view of (\ref{ee5}) and (\ref{ee4}) we obtain
\begin{eqnarray}\label{ee6}
&&R(X,Y)R(U,V)W-R(R(X,Y)U,V)W-R(U,R(X,Y)V)W\nonumber\\&&-R(U,V)R(X,Y)W=S(Y,R(U,V)W)X-S(X,R(U,V)W)Y\nonumber\\&&
-S(Y,U)R(X,V)W+S(X,U)R(Y,V)W-S(Y,V)R(U,X)W\nonumber\\&&+S(X,V)R(U,Y)W-S(Y,W)R(U,V)X+S(X,W)R(U,V)Y.
\end{eqnarray}
Substituting $X=U=\xi$ in (\ref{ee6}) and using (\ref{b7}),(\ref{b8}),(\ref{b9}) and (\ref{b10}) yields
\begin{eqnarray}\label{ee7}
&&-g(V,W)Y+g(Y,W)V-R(Y,V)W\nonumber\\&&=\eta(W)S(Y,V)\xi-2g(V,W)Y-2R(Y,V)W\nonumber\\&&
+2g(Y,W)\eta(V)\xi-S(Y,W)V
+S(Y,W)\eta(V)\xi\nonumber\\&&+2g(Y,V)\eta(W)\xi.
\end{eqnarray}
Taking the inner product of (\ref{ee7}) with $Z$ we obtain
\begin{eqnarray}\label{ee8}
&&g(R(Y,V)W,Z)+g(V,W)g(Y,Z)+g(Y,W)g(V,Z)\nonumber\\&&+S(Y,W)g(V,Z)-S(Y,V)\eta(W)\eta(Z)
-S(Y,W)\eta(V)\eta(Z)\nonumber\\&&-2g(Y,V)\eta(W)\eta(Z)-2g(Y,W)\eta(V)\eta(Z)=0.
\end{eqnarray}
Let $\{e_{i}\}(1\leq i\leq 3)$ be an orthonormal basis of the tangent space at any point. Now taking summation over $i=1,2,3$ of the relation (\ref{ee8}) for $V=W=e_{i}$ gives
\begin{eqnarray}\label{ee9}
S(Y,Z)=-2g(Y,Z).
\end{eqnarray}
Also, from (\ref{e2}) using Proposition 2.1 we infer
\begin{eqnarray}
  S(Y,Z)=-(\lambda+1) g(Y,Z)+(\lambda-1)\eta (Y)\eta (Z).\label{ee10}
  \end{eqnarray}
In view of  (\ref{ee9}) and (\ref{ee10}) we have
\begin{eqnarray}
  (\lambda-1) [g(Y,Z)-\eta (Y)\eta (Z)]=0.\label{ee11}
  \end{eqnarray}
It follows  $\lambda = 1$, $\mu = 1$, which leads to the following:
\begin{thm} Let $(M,\phi, \xi,\eta,g)$ be a Kenmotsu 3-manifold. If $(g,\xi,\lambda,\mu)$ is an $\eta$-Ricci soliton on $M$ and the curvature condition $R.R=Q(S,R)$ holds, then $\lambda =\mu = 1$ and the manifold is an Einstein manifold.\end{thm}
\section{Example of a proper $\eta$-Ricci Soliton on a Kenmotsu 3-manifold}
We consider the 3-dimensional manifold $M=\{(x,y,z)\in \mathbb{R}^{3}, z\neq0\},$ where
$(x,y,z)$ are the standard coordinate of $\mathbb{R}^{3}.$
The vector fields$$e_{1}=z\frac{\partial }{\partial
x},\hspace{7pt}e_{2}=z\frac{\partial } {\partial y}
,\hspace{7pt}e_{3}=-z\frac{\partial }{\partial z}$$ are linearly
independent at each point of $M.$
Let $g$ be the Riemannian metric defined by
$$g(e_{1},e_{3})=g(e_{1},e_{2})=g(e_{2},e_{3})=0,$$
$$g(e_{1},e_{1})=g(e_{2},e_{2})=g(e_{3},e_{3})=1.$$
Let $\eta $ be the 1-form defined by $\eta (Z)=g(Z,e_{3})$, for any
$Z\in \Gamma (TM).$
Let $\phi $ be the $(1,1)$-tensor field defined by $$\phi
e_{1}=-e_{2},\hspace{7pt} \phi e_{2}=e_{1},\hspace{7pt}\phi
e_{3}=0.$$
Then using the linearity of $\phi $ and $g$, we have $$\eta
(e_{3})=1,$$ $$\phi ^{2}Z=-Z+\eta (Z)e_{3},$$ $$g(\phi Z,\phi
W)=g(Z,W)-\eta (Z)\eta (W),$$ for any $Z,W\in \Gamma (TM).$\\

Then for $e_{3}=\xi $, $(\phi ,\xi ,\eta ,g)$
defines an almost contact metric structure on $M$.\\

Let $\nabla $ be the Levi-Civita connection with respect to the metric $g$. Then we have
\begin{eqnarray}
[e_{1},e_{3}]&=&e_{1}e_{3}-e_{3}e_{1}\nonumber\\&=&z\frac{\partial
}{\partial x}
(-z\frac{\partial }{\partial z})-(-z\frac{\partial }{\partial z})(z\frac{\partial }{\partial x})\nonumber\\
&=&-z^{2}\frac{\partial ^{2}}{\partial x\partial
z}+z^{2}\frac{\partial ^{2}}{\partial z\partial x}+
z\frac{\partial }{\partial x}\nonumber\\&=&e_{1}.\end{eqnarray} Similarly
$$[e_{1},e_{2}]=0\hspace{10pt}and\hspace{10pt}
[e_{2},e_{3}]=e_{2}.$$\\

The Riemannian connection $\nabla $ of the metric $g$ is given by
\begin{eqnarray} 2g(\nabla
_{X}Y,Z)&=&Xg(Y,Z)+Yg(Z,X)-Zg(X,Y)\nonumber\\&-&g(X,[Y,Z])-g(Y,[X,Z])
+g(Z,[X,Y]),\label{f1}\end{eqnarray} which is known as Koszul's formula.\\

Using (\ref{f1}) we have \begin{eqnarray} 2g(\nabla
_{e_{1}}e_{3},e_{1})&=&-2g(e_{1},-e_{1})\nonumber\\&=&
2g(e_{1},e_{1}).\label{f2}\end{eqnarray}

Again by (\ref{f1}) \begin{eqnarray} 2g(\nabla
_{e_{1}}e_{3},e_{2})=0=2g(e_{1},e_{2})\label{f3}\end{eqnarray}
and
\begin{eqnarray} 2g(\nabla _{e_{1}}e_{3},e_{3})=0=2g(e_{1},e_{3}).\label{f4}\end{eqnarray}

From (\ref{f2}), (\ref{f3}) and (\ref{f4}) we obtain
$$2g(\nabla _{e_{1}}e_{3},X)=2g(e_{1},X),$$ for all $X\in \Gamma (TM).$\\

Thus $$\nabla _{e_{1}}e_{3}=e_{1}.$$\\

Therefore, (\ref{f1}) further yields
$$\nabla _{e_{1}}e_{3}=e_{1},\hspace{10pt}\nabla _{e_{1}}e_{2}=0,\hspace{10pt}
\nabla _{e_{1}}e_{1}=-e_{3},$$
$$\nabla _{e_{2}}e_{3}=e_{2},\hspace{10pt}\nabla _{e_{2}}e_{2}=e_{3},\hspace{10pt}
\nabla _{e_{2}}e_{1}=0,$$ \begin{eqnarray}\nabla
_{e_{3}}e_{3}=0,\hspace{10pt}\nabla _{e_{3}}e_{2}=0,\hspace{10pt}
\nabla _{e_{3}}e_{1}=0.\label{f5}\end{eqnarray}
It follows $\nabla _{X}\xi=X-\eta (X)\xi$, for
$\xi=e_{3}$. Hence the manifold is a Kenmotsu manifold. It is known that
\begin{eqnarray} R(X,Y)Z=\nabla _{X}\nabla _{Y}Z-\nabla _{Y}\nabla
_{X}Z-\nabla _{[X,Y]}Z. \label{f6}\end{eqnarray}
With the help of the above results and using (\ref{f6}), it can be
easily verified that
$$R(e_{1},e_{2})e_{3}=0,\hspace{10pt}R(e_{2},e_{3})e_{3}=-e_{2},\hspace{10pt}
R(e_{1},e_{3})e_{3}=-e_{1},$$
$$R(e_{1},e_{2})e_{2}=-e_{1},\hspace{10pt}R(e_{2},e_{3})e_{2}=e_{3},\hspace{10pt}
R(e_{1},e_{3})e_{2}=0,$$
$$R(e_{1},e_{2})e_{1}=e_{2},\hspace{10pt}R(e_{2},e_{3})e_{1}=0,\hspace{10pt}
R(e_{1},e_{3})e_{1}=e_{3}.$$\\
From the above expressions of the curvature tensor $R$ we obtain
\begin{eqnarray}
S(e_{1},e_{1})&=&g(R(e_{1},e_{2})e_{2},e_{1})+g(R(e_{1},e_{3})e_{3},e_{1})
\nonumber\\&=&-2.\end{eqnarray}
Similarly, we have $$S(e_{2},e_{2})=S(e_{3},e_{3})=-2.$$\\

In this case from (\ref{b16}), for $\lambda=-1$ and $\mu=3$ the data $(g,\xi,\lambda, \mu)$ is an $\eta$-Ricci soliton on $(M,\phi,\xi,\eta, g)$.

\section*{Acknowledgement}

 The authors are thankful to the referee for his/her valuable suggestions and comments towards
 the improvement of the paper.


\begin{bibdiv}
\begin{biblist}


\bib{bla1}{book}{
title={Lecture notes in Mathematics}, author={Blair, D. E.,},
date={509,1976}, publisher={Springer-Verlag Berlin} }


\bib{bla2}{book}{
title={Riemannian Geometry of contact and symplectic manifolds}, author={Blair, D. E.,}, date={2002}, publisher={Birkh{\"a}user}, address={Boston} }

\bib{blaga1}{article}{title={$\eta$-Ricci solitons on Lorentzian para-Sasakian manifolds}, author={Blaga, A. M.,}, journal={Filomat}, volume={30}, date={2016},
pages={no. 2, 489-496} }

\bib{blaga2}{article}{title={$\eta$-Ricci solitons on para-Kenmotsu manifolds}, author={Blaga, A. M.,}, journal={Balkan J. Geom. Appl.}, volume={20}, date={2015},
pages={1-13} }

\bib{CC}{article}{title={From the Eisenhart problem to Ricci solitons in f-Kenmotsu manifolds}, author={Calin, C.,}, author={Crasmareanu, M.,}, journal={Bull. Malays. Math. Soc.}, volume={33(3)}, date={2010},
pages={361-368} }

\bib{chave1}{article}{title={Quasi-Einstein metrics and their renormalizability properties}, author={Chave, T.,}, author={Valent, G.,}, journal={Helv. Phys. Acta.}, volume={69}, date={1996},
pages={344-347} }

\bib{chave2}{article}{title={On a class of compact and non-compact quasi-Einstein metrics and their renormalizability properties}, author={Chave, T.,}, author={Valent, G.,}, journal={Nuclear Phys. B.}, volume={478}, date={1996},
pages={758-778} }

\bib{cho1}{article}{title={Ricci solitons and real hypersurfaces in a complex space form}, author={Cho, J. T.}, author={Kimura, M.,}, journal={Tohoku Math. J.}, volume={61, no. 2}, date={2009},
pages={205-212} }


\bib{de}{article}{title={On 3-dimensional Kenmotsu manifold}, author={De, U. C.,}, author={Pathak, G.,}, journal={ Bull. Malays. Math. Soc.}, volume={29}, date={2006},
pages={51-57} }

\bib{dd1}{article}{title={On semi-Riemannian manifolds satisfying the condition R.R=Q(S,R)}, author={Defever, F.,}, author={Deszcz, R.,}, journal={Geometry and topology of submanifolds, World Sci., Singapore}, volume={3}, date={1990},
pages={108-130} }

\bib{dd2}{article}{title={Some results on geodesic mappings of Riemannian manifolds satisfying the condition R.R=Q(S,R)}, author={Defever, F.,}, author={Deszcz, R.,}, journal={Periodica Mathematica Hungarica}, volume={29}, date={1994},
pages={267-276} }

\bib{dd3}{article}{title={ On conformally flat Riemannian manifolds satisfying certain curvature conditions}, author={Deszcz, R.,}, journal={Tensor (N.S.)}, volume={49}, date={1990}, pages={134-145} }

\bib{deshmukh1}{article}{title={Jacobi-type vector fields on Ricci solitons}, author={Deshmukh, S.,}, journal={Bull. Math. Soc. Sci. Math. Roumanie},
volume={55(103), 1}, date={2012},
pages={41-50} }

\bib{deshmukh2}{article}{title={A Note on Ricci Soliton}, author={Deshmukh, S.,}, author={Alodan, H.}, author={Al-Sodais, H.,}, journal={Balkan J.
                               Geom. Appl.}, volume={16, 1}, date={2011},
pages={48-55} }
\bib {dugg}  {article}{title={ Almost Ricci Solitons and Physical Applications},author={Duggal, K.L.,}, journal={International Electronic Journal of Geometry}, volume={2}, date={2017},
pages={1-10} }
\bib{friedam}{article}{title={  Non linear models in $2+\epsilon$ dimensions}, author={Friedan, D.,}, journal={Ann. Phys.}, volume={163}, date={1985},
pages={318-419} }

\bib{gray}{article}{title={Einstein-like manifolds which are not Einstein}, author={Gray, A.,}, journal={Geom. Dedicata}, volume={7}, date={1978},
pages={259-280} }

\bib{hr}{article}{title={The Ricci flow on surfaces, Mathematics and general relativity}, author={Hamilton, R. S.,}, journal={}, journal={(Contemp. Math. Santa Cruz, CA, 1986),  American Math. Soc.,} date={1988}, pages={237-263}}
\bib{ian} {article}{title={ Some remarkable structures on the products of an almost contact metric manifold with the real line}, author={Ianus, S.,},author={Smaranda, D.}, journal={ Papers from the National Coll. on Geometry and Topology, Univ. Timisoara}, date={1997}, pages={ 107-110}}
\bib{it}{article}{title={Ricci solitons on compact 3-manifolds}, author={Ivey, T.,}, journal={Diff. Geom. Appl.}, volume={3}, date={1993},
pages={301-307} }

\bib{ken}{article}{title={ A class of almost contact Riemannian manifolds}, author={Kenmotsu, K.,}, journal={Tohoku Math. J.}, volume={24}, date={1972},
pages={93-103} }

\bib{oub}{article}{title={New classes of almost contact metric structures}, author={Oubina, J.A.,}, journal={Publ. Math. Debrecen}, volume={32}, date={1985},
pages={187-193} }


\bib{dg}{article}{title={$\eta$-Ricci solitons on para-Sasakian manifolds}, author={Prakasha, D. G.,}, author={Hadimani, B. S.,}, journal={J. Geom.}, volume={108}, date={2017},
pages={383-392} }
\bib{tan} {article}{title={The automorphism groups of almost contact Riemannian manifolds}, author={Tanno, S.,}, journal={Tohoku Math. J.}, volume={21}, date=
{1969}, pages={221-228} }





\end{biblist}
\end{bibdiv}
\end{document}